\documentclass{amsart}
\usepackage{amsmath,amssymb,amsthm,amsfonts,mathrsfs,textcomp,mathtools,accents,color}
\usepackage{graphicx}
\usepackage{enumerate}
\usepackage{subcaption}

\makeatletter
\@namedef{subjclassname@2020}{%
  \textup{2020} Mathematics Subject Classification}
\makeatother

\usepackage{fullpage}

\newtheorem{theorem}{Theorem}[section]

\newtheorem{lemma}[theorem]{Lemma}
\newtheorem{proposition}[theorem]{Proposition}
\theoremstyle{definition}

\theoremstyle{remark}
\newtheorem{rem}[theorem]{Remark}
\theoremstyle{definition}

\numberwithin{equation}{section}

\newcommand{\set}[1]{\left\{#1\right\}}
\newcommand{\virg}[1]{``#1"}
\newcommand{\hyper}{_2F_1}
\newcommand{\R}{\mathbb R}
\newcommand{\Z}{\mathbb Z}
\newcommand{\N}{\mathbb N}
\newcommand{\Q}{\mathbb Q}
\newcommand{\C}{\mathbb C}
\newcommand{\Hp}{\mathbb H}

\newcommand{\PPP}{\mathcal{P}}
\newcommand{\QQQ}{\mathcal{Q}}

\newcommand{\MMM}{\mathcal{M}}
\newcommand{\BBB}{\mathcal{B}}

\begin{document}

\title{On the generalised transfer operators of the Farey map with complex temperature}
\author{Claudio Bonanno}
\address{Dipartimento di Matematica , Universit\`a di
Pisa, Largo Bruno Pontecorvo 5, 56127 Pisa, Italy} 
\email{claudio.bonanno@unipi.it} 

\date{}

\thanks{The author is partially supported by 
the PRIN Grant 2017S35EHN ``Regular and stochastic behaviour in dynamical 
systems'' of the Italian Ministry of University 
and Research (MUR), Italy. This research is part of the author's activity within 
the UMI Group ``DinAmicI'' (\texttt{www.dinamici.org}) and the 
Gruppo Nazionale di Fisica Matematica, INdAM, Italy.}

\begin{abstract}
We consider the problem of showing that 1 is an eigenvalue for a family of generalised transfer operators of the Farey map. This problem is related to the spectral theory of the modular surface via the Selberg Zeta function and the theory of dynamical zeta functions of maps. After briefly recalling these connections, we show that the problem can be formulated for operators on an appropriate Hilbert space and translated into a linear algebra problem for infinite matrices.
\end{abstract}
\maketitle

\section{Introduction} \label{intro}
In the modern theory of dynamical systems one of the most successful methods for the study of the statistical properties of a system is the \emph{thermodynamic formalism} approach. Originated from the works of Sinai, Ruelle, and Bowen, this approach has led to important results especially for uniformly hyperbolic systems, and has been largely developed in the last thirty years. We refer to \cite{bowen,ruelle-book}, classical books with recent editions, for the basic description of this method.

An important tool in the thermodynamic formalism approach is the \emph{transfer operator} of a dynamical system, a linear operator describing the evolution of measures on the phase space under the action of the system. In this paper we consider the transfer operator of the \emph{Farey map}, a map on the unit interval which is related to the continued fraction expansion of real numbers. The Farey map is interesting also from the dynamical systems point of view, being a prototypical example of a non-uniformly hyperbolic map of the interval and having the property that the only invariant measure $\mu$ which is absolutely continuous with respect to the Lebesgue measure is infinite, that is $\mu(0,1)=\infty$. These characteristics imply that the classical methods to study the statistical properties of the map fail, as for example the approach using the spectral properties of the transfer operator. Nevertheless more recent methods have been introduced that show that the transfer operator is useful to obtain the properties of a dynamical system also in these less-standard situations (see e.g. \cite{mel-ter,BGL} and the review \cite{gouezel}).

In this paper we first give a review of the results obtained by the author with collaborators on the transfer operator of the Farey map.  In fact, we consider a family of transfer operators which depend on a real parameter, sometimes called \emph{temperature} in the thermodynamic formalism approach.  In particular, we study the classification of the eigenfunctions of these operators in a Banach space of holomorphic functions. The main idea is that it is possible to read the action of these transfer operators on a Hilbert space of functions defined on $\R^+$ for which the (generalised) Laguerre polynomials are an orthogonal basis. This has led to the possibility of studying the spectral properties of these operators via a matrix approach (see \cite{isola,BGI,imen1,imen2}). Then we extend this approach to the case of complex temperatures, an important case to consider in light of the connections with the spectral theory of the hyperbolic surfaces developed by Mayer and Lewis-Zagier (see \cite{Ma2,lewis,CM,lewis-zagier,BLZ2,BLZ1} and the expository paper \cite{PZ} which contains a description of the extensions of these connections to other quotients of the hyperbolic plane). A side result of these connections is that one can formulate a functional equation whose solutions give information on the position of the non-trivial zeroes of the Riemann Zeta function. Thus we obtain yet another formulation of the Riemann hypothesis in this case as a linear algebra problem for infinite matrices.

The paper is structured as follows. In Section \ref{sec:setting} we introduce the Farey map and the family of generalised transfer operators we consider. Then we describe the relations between the Farey map and the continued fraction expansion of real numbers and with the spectral theory of the hyperbolic surfaces. In Section \ref{sec:result} we state the results on the classification of the eigenfunctions of the transfer operators which lead to the formulation of the eigenvalue problem on the appropriate Hilbert space of functions. In particular, it turns out that it interesting to find the values of the temperature for which the transfer operator has 1 as eigenvalue with eigenfunction of a given type. Finally we describe the matrix approach in the case of complex temperatures. The main results are the formulation of the problems in terms of infinite matrices described in Theorems \ref{thm:mat-0} and \ref{thm:mat-b}. 

\section{The generalised transfer operators of the Farey map} \label{sec:setting}
Let $J_0:=[0,1/2]$ and $J_1:=[1/2,1]$, the \emph{Farey map} is the continuous transformation of the unit interval $F:[0,1]\to [0,1]$ defined to be
\begin{equation} \label{eq:farey}
F(x):= \left\{ \begin{array}{ll} \frac{x}{1-x}, & \text{if $x\in J_0$;}\\[0,2cm] \frac{1-x}{x}, & \text{if $x\in J_1$;}
\end{array}
\right.
\end{equation} 
which is piecewise monotone and differentiable, and with full branches (i.e. $F(J_0)=F(J_1)=[0,1]$). Let $\phi_0$ and $\phi_1$ denote the two local inverses of $F$, that is
\[
\begin{aligned}
& \phi_0 := \left( F|_{J_0} \right)^{-1} : [0,1] \to \left[ 0,\frac 12\right],\qquad \phi_0(x) = \frac{x}{1+x},\\[0.3cm]
& \phi_1 := \left( F|_{J_1} \right)^{-1} : [0,1] \to \left[ \frac 12,1\right],\qquad \phi_1(x) = \frac{1}{1+x}.
\end{aligned}
\]
The \emph{transfer operator} of the Farey map is the positive linear operator $\PPP$ which acts on the space of measurable functions on the unit interval by
\begin{equation}\label{eq:to}
(\PPP f)(x) := |\phi'_0(x)|\, f(\phi_0(x)) +  |\phi'_1(x)|\, f(\phi_1(x)) = \frac{1}{(1+x)^2}\, f\left(\frac{x}{1+x}\right) + \frac{1}{(1+x)^2}\, f\left(\frac{1}{1+x}\right).
\end{equation}
The main feature of the transfer operator of a map is its relation with the absolutely continuous measures on $[0,1]$ which are invariant for the action of the map. For the Farey map it is known that the only absolutely continuous invariant measure $\mu$ has density $h(x) = 1/x$ and is not normalisable, i.e. $\mu(0,1)=\infty$. It is straightforward to verify that $\PPP h = h$. In addition one can show that $h$ is the only function fixed by $\PPP$. Other spectral properties of the transfer operator of a map are related to the statistical properties of the associated dynamical system, in particular to the property of mixing and the speed of decay of correlations (see e.g. \cite{mel-ter,gouezel,BGL}).

In this paper we are interested in the family of \emph{signed generalised transfer operators} of the Farey map. For $q\in \C$, let us consider the generalised version of the two operators appearing in \eqref{eq:to} separately,
\[
\begin{aligned}
& (\PPP_{0,q} f)(x) := |\phi'_0(x)|^q\, f(\phi_0(x)) = \frac{1}{(1+x)^{2q}}\, f\left(\frac{x}{1+x}\right), \\[0.2cm]
& (\PPP_{1,q} f)(x) := |\phi'_1(x)|^q\, f(\phi_0(x)) = \frac{1}{(1+x)^{2q}}\, f\left(\frac{1}{1+x}\right),
\end{aligned}
\]
that is, we consider separately the contributions coming from the two branches of the Farey map. The two operators have different spectral properties when acting on a suitable Banach space of functions. This is related to the different dynamical behaviour of the two branches of the map, the \emph{parabolic} one, $F$ restricted to $J_0$, and the \emph{expanding} one, $F$ restricted to $J_1$.

The signed generalised transfer operators of the Farey map are then defined to be the operators acting on the space of measurable functions by
\begin{equation}\label{eq:toq}
(\PPP_q^\pm f)(x) := (\PPP_{0,q} f)(x) \pm (\PPP_{1,q} f)(x), \quad q\in \C.
\end{equation}
The study of generalised transfer operators with real parameter $q$ is standard in the thermodynamic formalism approach to dynamical systems, see \cite{ruelle-book}, in which $q$ is called \emph{temperature}. The motivations for studying the signed operators with complex temperature $q$ come from the connections between the continued fraction expansion of real numbers and the spectral theory of the modular surface which we briefly recall now.

\subsection{The Farey map, continued fractions and the modular surface} \label{subsec:selberg}

Given a real number $x$, its \emph{(regular) continued fraction expansion} is obtained by writing
\[
x = [a_0; a_1, a_2, \dots] := a_0 + \frac{1}{a_1 + \frac{1}{a_2 + \frac{1}{\ddots}}}
\]
with $a_0 := \lfloor x \rfloor \in \Z$, and $a_n\in \N$ for all $n\ge 1$. It is well-known that the expansion is finite if and only if $x\in \Q$, and it is convergent for $x\in \R\setminus \Q$. A classical reference where to find more details is \cite{cf}. Here we are interested in the dynamical systems approach. Let $G:[0,1]\to [0,1]$ be the transformation of the unit interval defined by
\begin{equation} \label{eq:gauss}
G(x) := \left\{ \begin{array}{ll} \frac 1x - \left\lfloor \frac 1x \right\rfloor, & \text{if $x\in (0,1]$;}\\[0.2cm] 0, & \text{if $x=0$.} \end{array} \right.
\end{equation}
The map $G$ is known as the \emph{Gauss map} and is strictly related to the continued fraction expansion of a real number. It is immediate to show that
\[
[0,1]\ni x= [0; a_1, a_2, a_3,\dots] \mapsto G(x) = [0; a_2, a_3,\dots],
\]
so that the Gauss map acts as a shift map on the sequence of the continued fraction coefficients $\{a_n\}$ of $x$, and for all $x\in [0,1]$
\[
a_n = \left\lfloor \frac{1}{G^{n-1}(x)} \right\rfloor, \quad \forall \, n\ge 1,
\]
where as usual $G^n := G\circ G\circ \dots \circ G$ for $n$ times, and $G^0=Id$.

The Gauss map is piecewise differentiable and invertible with respect to the countable partition $\{A_n\}$, $n\ge 1$, of the unit interval given by
\[
A_n := \left( \frac{1}{n+1}, \frac 1n \right].
\]
If we set $\psi_n$ to denote the local inverse of $G$, that is
\[
\psi_n := \left( G|_{A_n} \right)^{-1} : [0,1) \to A_n,\qquad \psi_n(x) = \frac{1}{n+x},
\]
the generalised transfer operators of the Gauss map are the family of linear operators $\QQQ_q$, with $q\in \C$, acting on the space of measurable functions on the unit interval by
\[
(\QQQ_q g)(x) := \sum_{n=1}^\infty\, |\psi'_n(x)|^q\, g(\psi_n(x)) = \sum_{n=1}^\infty\, \frac{1}{(n+x)^{2q}}\, g\left(\frac{1}{n+x}\right).
\]
A first study of these operators can be found in \cite{Ma1}, but it is \cite{Ma2} the paper which showed the beautiful connection of these operators with the spectral theory of the modular surface.

Let $\Hp=\{ z=x+iy\, :\, y>0\}$ be the Poincar\'e half-plane, that is the upper half-plane endowed with the hyperbolic metric which makes it a surface with constant negative curvature. The group of order preserving isometries of the Poincar\'e half-plane is given by the quotient $PSL(2,\R) := SL(2,\R)/\pm I$ of matrices $M\in SL(2,\R)$ which act on $\Hp$ as M\"obius transformations
\[
\Hp \ni z \mapsto M(z) = \begin{pmatrix} a & b \\[0.1cm] c & d \end{pmatrix}\, \cdot\, z = \frac{az+b}{cz+d} \in \Hp,
\]
from which it is evident that $M$ and $-M$ represent the same transformation. An important discrete subgroup of $PSL(2,\R)$ is the \emph{modular group} $PSL(2,\Z)$ of matrices with integer entries. The \emph{modular surface} $\MMM:= PSL(2,\Z)\backslash \Hp$ is the quotient of $\Hp$ by left actions of matrices in $PSL(2,\Z)$. It can be seen as a non-compact finite measure surface with constant negative curvature, with a cusp end and two conical singularities. Let 
\[
\triangle_\Hp := -y^2\, \left( \frac{\partial^2}{\partial x^2} +  \frac{\partial^2}{\partial y^2} \right)
\]
be the Laplace-Beltrami operator on $\Hp$, a classical problem is the study of the spectrum of $\triangle_\Hp$ on $\MMM$, or equivalently the spectrum of $\triangle_\Hp$ restricted to $PSL(2,\Z)$-invariant functions on $\Hp$. Without entering into a complete description of the spectral resolution of $L^2(\MMM)$ (see e.g. \cite{terras} for more details), we recall that $\triangle_\Hp$ restricted to $L^2(\MMM)$ has a continuous spectrum given by the line $[1/4, +\infty)$, and discrete spectrum given by $0$ and a sequence of values $\{ \lambda_n\}$ embedded in the continuous spectrum. 

The results in \cite{Ma2} use the theory of the Selberg Zeta function and of the dynamical zeta function of the Gauss map. Recalling all these details is beyond the aims of this paper, so we limit ourselves to state a corollary of the main result in \cite{Ma2} for the spectral theory of $\MMM$, as sharpened after \cite{efrat}. Let
\[
D := \set{ z\in \C\, :\, |z-1| < \frac 32}
\]
and denote by $A_\infty(D)$ the Banach space of holomorphic functions in $D$ which can be extended to continuous functions on $\bar{D}$ equipped with the supremum norm. There is a relation between the spectrum of $\triangle_\Hp$ on $\MMM$ and the functions in $A_\infty(D)$ which are eigenfunctions of the generalised transfer operator $\QQQ_q$ of the Gauss map with eigenvalues $\pm 1$.

\begin{theorem}[\cite{Ma2},\cite{efrat}] \label{thm:mayer}
Let $q\in \C$ with $\Re(q)>0$ and $q\not=1/2$. Then:
\begin{enumerate}[(i)]
\item there exists a nonzero $g\in A_\infty(D)$ such that $\QQQ_{q}g = g$ if and only if either $\lambda_q := q(1-q)$ is in the discrete spectrum of $\triangle_\Hp$ restricted to $L^2(\MMM)$ with eigenfunction $u_{\lambda_q}$ satisfying $u_{\lambda_q}(-x+iy) = u_{\lambda_q}(x+iy)$, or $2q$ is a non-trivial zero of the Riemann Zeta function;
\item there exists a nonzero $g\in A_\infty(D)$ such that $\QQQ_{q}g = - g$ if and only if $\lambda_q := q(1-q)$ is in the discrete spectrum of $\triangle_\Hp$ restricted to $L^2(\MMM)$ with eigenfunction $u_{\lambda_q}$ satisfying $u_{\lambda_q}(-x+iy) = - u_{\lambda_q}(x+iy)$.
\end{enumerate}
\end{theorem}

The values of $q$ for which $2q$ is a non-trivial zero of the Riemann Zeta functions have a spectral and a dynamical interpretation in terms of the scattering matrix of the modular surface (see e.g. \cite{andreas}). In addition, since the eigenvalues $\lambda_q$ are real and consist of 0 and of a sequence embedded in $[1/4,+\infty)$, it follows that the values of $q$ for which there exists a nonzero $g\in A_\infty(D)$ such that $\QQQ_{q}g = \pm g$ are $1$ or satisfy $\Re(q) = 1/2$. This is well known also from the work of Selberg (see e.g. \cite{iwaniec}).

Let us now introduce the Farey map in this framework. The idea is to notice that the Gauss map \eqref{eq:gauss} may be defined by using the Farey map \eqref{eq:farey} and inducing, a classical idea of dynamical systems theory. In many situations, it is possible to study the dynamical behaviour of a system by looking at the passages of the orbits through a subset of the phase space. This is for example the case when the dynamics has regular properties on the whole phase space except for a positive measure set, on which the dynamics is chaotic. Looking at the Farey map, we have already observed that the dynamics is regular when restricted to $J_0=[0,1/2]$, and the chaotic behaviour of the map is determined by the dynamics on $J_1=[1/2,1]$. Then one may want to ``accelerate'' the dynamics, by looking at the orbits only after a passage through $J_1$. This leads to the notion of jump transformation.

Let
\[
\tau : [0,1] \to \N,\quad \tau(x) := 1 + \min \set{ k\ge 0\, :\, F^k(x) \in \left(\frac 12,1\right]},
\]
with the convention that $\tau(x) = \infty$ if the orbit of $x$ does not visit $(1/2,1]$. Then the Gauss map $G$ is nothing but the jump transformation of the Farey map $F$ on the set $(1/2,1]$ given by $G(x) = F^{\tau(x)}(x)$ if $\tau(x)<\infty$, and $G(x)=0$ otherwise. Note that the level sets of the function $\tau$ are exactly the sets $A_n=(1/(n+1),1/n]$, that is $\tau(x)=n$ if and only if $x\in A_n$. 

The definition of $G$ as jump transformation of $F$ leads to the following relations among the two maps and among their local inverses:
\[
G|_{A_n} = F|_{J_1}\circ (F|_{J_0})^{n-1},\quad \psi_n =  \left( G|_{A_n} \right)^{-1} = \phi_0^{n-1} \circ \phi_1,\quad \forall\, n\ge 1.
\]
A formal relation between the generalised transfer operators of the two maps follows. For $n=1$ we have
\[
|\psi_1'(x)|^q\, g(\psi_1(x)) = |\phi_0'(x)|^q\, g(\phi_0(x)) = (\PPP_{0,q} g)(x),
\]
for $n=2$ the operator $\PPP_{1,q}$ appears in the formula,
\[
|\psi_2'(x)|^q\, g(\psi_2(x)) = |\phi_0'(\phi_1(x))|^q\, |\phi_1'(x)|^q\, g(\phi_0(\phi_1(x))) = (\PPP_{1,q} (\PPP_{0,q} g))(x),
\]
and for all $n\ge 1$ one can show by induction that
\[
|\psi_n'(x)|^q\, g(\psi_n(x)) = (\PPP_{1,q} (\PPP_{0,q}^{n-1} g))(x).
\]
Hence
\[
\QQQ_q g = \sum_{n=1}^\infty\, \PPP_{1,q} (\PPP_{0,q}^{n-1} g)
\] 
as a formal expression, whose convergence depends on the properties of the function $g$. We can then write that
\begin{equation}\label{relation-fg-to-prel}
(1\mp \QQQ_q) (1-\PPP_{0,q}) = 1-\PPP_q^\pm.
\end{equation}
In the next section we recall from \cite{BI1} the results on the eigenfunctions of the operators $\PPP_q^\pm$ which are used to formulate \eqref{relation-fg-to-prel} on the right functional spaces. We will then be able to state the analog of Theorem \ref{thm:mayer} connecting the spectral properties of $\PPP_q^\pm$ and of the modular surface $\MMM$. 

\section{The eigenvalue-1 problem} \label{sec:result}

The first half of \cite{BI1} is devoted to the characterisation of the analytic eigenfunctions of the operators $\PPP_q^\pm$ for $q\in \C$ with $\Re(q)>0$ and $q\not=1/2$. Here we recall a simplified version of these results, considering only the case of the eigenvalue $1$. As discussed before, the unique absolutely continuous invariant measure of the Farey map has density $h(x)=1/x$, and it satisfies $\PPP_1^+h=h$. Therefore it is natural to consider the space $H(B)$ of holomorphic functions on the domain
\[
B:= \set{ z\in \C\, :\, |z-1| < \frac 12}
\]
and to look for solutions of $\PPP_q^\pm f=f$ on $H(B)$. 

Let's introduce some notations. First, for $a\in \C$ we introduce the symbol $\chi_a$ to denote the function
\[
\chi_a : \R^+ \to \C,\quad \chi_a(t) := t^a,
\]
and $m_a$ to denote the absolutely continuous measure on $\R^+$ with density $dm_a(t)/dt = \chi_{2a-1}(t) e^{-t}$. For $q\in \C$ we let $\xi:= \Re(q)>0$, and for all $p\in [1,\infty]$ consider the Banach spaces $(L^p(m_\xi), \|\cdot\|_p)$ defined as
\begin{equation} \label{mq}
L^p(m_\xi) := \set{\phi :\R^+\to \C\, :\, \|\phi\|_p := \int_0^\infty\, |\phi(t)|^p\, dm_\xi(t) <\infty}.
\end{equation}
Notice that $L^p(m_\xi) \subseteq L^1(m_\xi)$ for all $p\in [1,\infty]$. Given $\phi\in L^1(m_\xi)$ we recall the integral transform $\BBB_q$ introduced in \cite{isola} defined by
\[
\BBB_q[\phi](z):= \chi_{-2q}(z)\, \int_0^\infty\, e^{-\frac tz}\, \phi(t)\, dm_q(t),
\]
which is continuous as an operator from $(L^1(m_\xi), \|\cdot\|_1)$ to $H(B)$ with the topology induced by the family of supremum norms on the compact subsets of $B$. We also recall that the application of the $\BBB_q$-transform may be extended for $q\not=1/2$ to the function $\chi_{-1}(t)$ as in \cite[Rem. 2.4]{BI1}, and one obtains
\[
\BBB_q\left[\chi_{-1} \right](z) = \Gamma(2q-1)\, \chi_{-1}(z),
\]
where $\Gamma(\cdot)$ denotes the Gamma function.

We can now summarize the results of \cite{BI1} on the analytic solutions of $\PPP_q^\pm f=f$.
\begin{theorem}[\cite{BI1}] \label{thm:bi}
Let $q\in \C$ with $\xi=\Re(q)>0$ and $q\not=1/2$. If $f\in H(B)$ and $\PPP_q^\pm f = f$, then $f(z)$ is holomorphic on the half-plane $\{\Re(z)>0\}$ and there exists $\phi \in L^2(m_\xi)$ with $\phi(0)$ finite and $\phi(t) = \phi(0) + O(t)$ as $t\to 0^+$, such that
\begin{equation} \label{forma-giusta}
f(z) = \frac{c}{z^{2q}} + \BBB_q\left[ b\,\chi_{-1}+\phi\right](z),
\end{equation}
for $c,b \in \C$. If $f$ solves $\PPP_q^-f=f$ then $b=0$.
\end{theorem}
 
Equation \eqref{relation-fg-to-prel} turns out to be correct when written for functions written as in \eqref{forma-giusta}. This is the first step in \cite{BI1} to write the result analogous to Theorem \ref{thm:mayer}. Another idea is to use the equivalence between the equation $\PPP_q^\pm f=f$ and the \emph{three term functional equation} studied in \cite{lewis,lewis-zagier} (see \cite{BI2} for an extension of this equivalence). The following result is proved in \cite{BI1}.  

\begin{theorem}[\cite{BI1}] \label{thm-boniso}
Let $q\in \C$ with $\xi=\Re(q)>0$ and $q\not=1/2$. If $f(z)$ is a function of the form \eqref{forma-giusta} with $\phi\in L^2(m_\xi)$ and $c,b\in \C$, then
\[
(1\mp \QQQ_q)(1-\PPP_{0,q}) f = (1-\PPP_q^\pm) f \pm c.
\]
In particular, the operator $\QQQ_q$ has an eigenfunction $g\in A_\infty(D)$ with eigenvalue $\lambda_\QQQ =\pm 1$ if and only if $\PPP_q^\pm$ has an eigenfunction $f$ written as in \eqref{forma-giusta} with $\phi \in L^2(m_\xi)$ and $c=0$. It follows that:
\begin{enumerate}[(i)]
\item there exists a nonzero $f(z)$ written as in \eqref{forma-giusta} with $\phi\in L^2(m_\xi)$ and $c=b=0$ such that $\PPP_{q}^+f = f$ if and only if $\lambda_q := q(1-q)$ is in the discrete spectrum of $\triangle_\Hp$ restricted to $L^2(\MMM)$ with eigenfunction $u_{\lambda_q}$ satisfying $u_{\lambda_q}(-x+iy) = u_{\lambda_q}(x+iy)$;
\item there exists a nonzero $f(z)$ written as in \eqref{forma-giusta} with $\phi\in L^2(m_\xi)$ and $c=b=0$ such that $\PPP_{q}^-f = f$ if and only if $\lambda_q := q(1-q)$ is in the discrete spectrum of $\triangle_\Hp$ restricted to $L^2(\MMM)$ with eigenfunction $u_{\lambda_q}$ satisfying $u_{\lambda_q}(-x+iy) =-u_{\lambda_q}(x+iy)$;
\item there exists a nonzero $f(z)$ written as in \eqref{forma-giusta} with $\phi\in L^2(m_\xi)$, $c=0$ and $b\not=0$ such that $\PPP_{q}^+f = f$ if and only if $2q$ is a non-trivial zero of the Riemann Zeta function or $q=1$.
\end{enumerate}
\end{theorem}

The spectral properties of the Laplace-Beltrami operator on the modular surface $\MMM$ may then be studied by using the generalised transfer operator $\PPP_q^\pm$ of the Farey map. In particular, the eigenfunction equation for $\PPP_q^\pm$ is a functional equation with three terms, which is much easier to treat than the equivalent one for the operator of the Gauss map $\QQQ_q$ which has an infinite number of terms. Notice that the case $q=1$ corresponds to the transfer operator defined in \eqref{eq:to}, which has eigenfunction $f(z)=1/z$ with eigenvalue 1. Hence in this case $f$ is written as in \eqref{forma-giusta} with $\phi\equiv 0$ and $c=0$.

One possible way to study the existence of solutions to the equation $\PPP_q^\pm f=f$ in the correct functional space has been proposed in \cite{BGI,imen1} for the case of real values of $q$. In this paper we extend this approach to the general case of $q\in \C$ with $\xi=\Re(q)>0$ and $q\not=1/2$, which is the situation one needs to consider to work in the framework of Theorem \ref{thm-boniso}. 

Let us consider the linear operators $M$ and $N_q$ defined for $q\in \C$ with $\xi=\Re(q)>0$ on functions $\psi:\R^+\to \C$ by
\begin{equation}\label{eq-m-nq}
M(\psi)(t) := e^{-t}\, \psi(t),\qquad N_q(\psi)(t) := \int_0^\infty\, \frac{J_{2q-1}(2\sqrt{st})}{(st)^{q-1/2}}\, \psi(s)\, dm_q(s).
\end{equation}
Here $J_\nu(\cdot)$ denotes the \emph{Bessel function of first kind}, which has the power series expansion
\begin{equation}\label{bessel-exp}
J_\nu(t) = \frac{t^\nu}{2^\nu}\, \sum_{m=0}^\infty\, \frac{(-1)^m\, t^m}{2^{2m}\, m!\, \Gamma(m+\nu+1)}
\end{equation}
for $\Re(\nu)>-1$, and satisfies $J_\nu(t) = O(t^\nu)$ as $t\to 0^+$, and $J_\nu(t) = O(t^{-\frac 12})$ as $t\to \infty$ (see \cite[vol. II]{E}).

\begin{proposition}[\cite{BI1}] \label{prop-m-nq}
Let $q\in \C$ with $\xi=\Re(q)>0$ and $q\not=1/2$. The operators $M$ and $N_q$ are bounded operators from $L^2(m_\xi)$ to $L^2(m_\xi)$, and for all $\phi\in L^2(m_\xi)$ it holds
\[
\begin{aligned}
& \PPP_{0,q} \left( \BBB_q\left[ \chi_{-1}+\phi\right] \right) = \BBB_q\left[ M \left(\chi_{-1}+\phi\right) \right], \\
& \PPP_{1,q} \left( \BBB_q\left[ \chi_{-1}+\phi\right] \right) = \BBB_q\left[ N_q \left(\chi_{-1}+\phi\right) \right]. 
\end{aligned}
\]
In addition $N_q(\chi_{-1}) \in L^2(m_\xi)$.
\end{proposition}

We can then state points (i)-(iii) of Theorem \ref{thm-boniso} by using the operators $M$ and $N_q$. Hence we are interested in looking for nonzero functions $\phi\in L^2(m_\xi)$ which are solutions to 
\begin{equation}\label{the-problem-0}
( M\pm N_q - I) \phi = 0, 
\end{equation}
or to
\begin{equation}\label{the-problem-b}
( M+ N_q - I) \phi = (I-M-N_q) \chi_{-1}. 
\end{equation}
Solutions to \eqref{the-problem-0} give values $q$ corresponding to the point spectrum of the Laplace-Beltrami operator on the modular surface, and solutions to \eqref{the-problem-b} give values $q$ corresponding to the non-trivial zeroes of the Riemann Zeta function.  

\begin{rem}\label{rem-q1}
Using the expansion in \eqref{bessel-exp}, one shows that 
\[
N_q (\chi_{-1})(t) = \sum_{m=0}^\infty\, \frac{(-1)^m}{m!\, (m+2q-1)}\, t^m,
\]
see \cite[p. 903]{BI1}. Hence for $q=1$ one gets $N_1(\chi_{-1}) (t) = (1-e^{-t})\, \chi_{-1}(t)$, that is, $N_1(\chi_{-1}) = (I-M) \chi_{-1} $, from which it follows that \eqref{the-problem-b} for $q=1$ admits the solution $\phi\equiv 0$.
\end{rem}

\subsection{A matrix approach} \label{sec:matrix} 
In this section we formulate equations \eqref{the-problem-0} and \eqref{the-problem-b} in terms of infinite matrices. The main idea is to consider the space $L^2(m_\xi)$ in \eqref{mq} as a Hilbert space endowed with the scalar product
\[
(f,g)_\xi := \int_0^\infty\, f(t)\, \overline{g(t)}\, dm_\xi(t).
\]
As above we assume $\xi:=\Re(q)>0$ and $q\not= 1/2$. It turns out that it is useful to consider the basis of orthogonal polynomials on $L^2(m_\xi)$ given by the \emph{(generalised) Laguerre polynomials}
\begin{equation}\label{laguerre}
L^{2\xi-1}_n(t) := \sum_{\ell=0}^n\, \frac{\Gamma(n+2\xi)}{(n-\ell)!\, \Gamma(\ell+2\xi)}\, \frac{(-1)^\ell}{\ell!}\, t^\ell,\quad \forall\, n\in \N_0,
\end{equation}
which satisfy the relation
\[
(L^{2\xi-1}_n\, ,\, L^{2\xi-1}_m)_\xi = \frac{\Gamma(n+2\xi)}{n!}\, \delta_{n,m}, \quad \forall\, n,m\in \N_0,
\]
where $\delta_{\cdot,\cdot}$ is the Kronecker delta (see e.g. \cite{pol}) and $\Gamma(\cdot)$ denotes the Gamma function. Then, given a function $\phi\in L^2(m_\xi)$ we can write it as
\begin{equation}\label{expansion}
\phi(t) = \sum_{n=0}^\infty\, \Phi_n\, L_n^{2\xi-1}(t),
\end{equation}
where $\{\Phi_n\}_n$ is a sequence in $\C$ such that
\[
(\phi\, ,\, L_n^{2\xi-1})_\xi = \Phi_n \frac{\Gamma(n+2\xi)}{n!} \quad \text{and} \quad \|\phi\|_2^2 =  \sum_{n=0}^\infty\, |\Phi_n|^2\, \frac{\Gamma(n+2\xi)}{n!} <+\infty.
\]
We start considering problem \eqref{the-problem-0}. Using the basis of Laguerre polynomials, a function $\phi\in L^2(m_\xi)$ is a nonzero solution of 
\eqref{the-problem-0} if and only if
\[
\Big( (M\pm N_q)\phi\, ,\, L^{2\xi-1}_k\Big)_\xi = (\phi\, ,\, L_k^{2\xi-1})_\xi, \quad \forall\, k\in \N_0.
\]
Since $M$ and $N_q$ are linear operators, using \eqref{expansion},
\[
(M\pm N_q)\phi = \sum_{n=0}^\infty\, \Phi_n\, (M\pm N_q) L^{2\xi-1}_n,
\]
hence there exists a nonzero solution to \eqref{the-problem-0} if and only if there exists a nonzero infinite vector $\Phi := (\Phi_n)$ such that
\begin{equation}\label{eq-prob-0}
\left\{ \begin{array}{l}
\sum_{n=0}^\infty\, |\Phi_n|^2\, \frac{\Gamma(n+2\xi)}{n!}<\infty \\[0.2cm]
A^\pm_q \Phi = D_q \Phi
\end{array} \right.
\end{equation}
where $A^\pm_q$ is the infinite matrix defined by
\begin{equation}\label{inf-mat-0-impl}
A^\pm_q = (a^\pm_{kn}(q))_{k,n\in \N_0}, \quad a^\pm_{kn}(q) := \Big( (M\pm N_q)L^{2\xi-1}_n\, ,\, L^{2\xi-1}_k\Big)_\xi,
\end{equation}
and $D_q$ is the infinite diagonal matrix given by
\begin{equation}\label{inf-mat-diag}
D_q = (d_{kn}(q))_{k,n\in \N_0}, \quad d_{kn}(q) := \frac{\Gamma(k+2\xi)}{k!}\, \delta_{k,n}.
\end{equation}
We now evaluate the entries $a^\pm_{kn}(q)$ of the matrix $A^\pm_q$. The part concerning the operator $M$ has been already considered in \cite{BGI,imen1}. Since it does not involve $\Im(q)$, the computation is the same as for $q$ real.

\begin{lemma}[\cite{BGI,imen1}] \label{lem:m-matr}
For all $\xi>0$ it holds
\[
\Big( M(L_n^{2\xi-1})\, ,\, L_k^{2\xi-1} \Big)_\xi = \frac{\Gamma(k+n+2\xi)}{k!\, n!}\, 2^{-(k+n+2\xi)},
\]
for all $k,n\in \N_0$.
\end{lemma}

The computations involving $N_q$ cannot be reduced to the case $q$ real. Let us recall the definition of the \emph{hypergeometric function} $\hyper\left( \alpha,\beta;\gamma;z \right)$ defined for $|z|<1$ and complex numbers $\alpha, \beta$ and $\gamma$ by
\[
\hyper\left( \alpha,\beta;\gamma;z \right) := \sum_{i=0}^\infty\, \frac{\Gamma(\alpha+i)\Gamma(\beta+i)\Gamma(\gamma)}{\Gamma(\alpha)\Gamma(\beta)\Gamma(\gamma+i)}\, \frac{z^i}{i!},
\]
(see e.g. \cite[vol. I]{E}).

\begin{lemma} \label{lem:n-matr}
Assume $\xi:=\Re(q)>0$ and $q\not= 1/2$. Then
\[
\begin{aligned}
& \Big( N_q(L_n^{2\xi-1})\, ,\, L_k^{2\xi-1} \Big)_\xi  =  \frac{\Gamma(n+2\xi)\Gamma(k+2\xi)}{2^{k+2\xi}}\, \sum_{\ell=0}^n\,  \frac{(-1)^\ell\, \Gamma(\ell+2q)}{(n-\ell)!\, \Gamma(\ell+2\xi)}\, \cdot \\
& \cdot \left(\sum_{j=0}^{\min\{\ell,k\}}\, \frac{2^{-j}}{j!\, (\ell-j)! \, (k-j)! \, \Gamma(j+2q)}\, \hyper\left( -\ell+j,2\xi +j;2q+j;\frac 12\right)\right),
\end{aligned}
\]
for all $k,n\in \N_0$. 
\end{lemma}

\begin{proof}
We start using \eqref{eq-m-nq} and \eqref{laguerre} to write
\[
N_q(L_n^{2\xi-1}) (t) = \sum_{\ell=0}^n\, \frac{(-1)^\ell\, \Gamma(n+2\xi)}{(n-\ell)!\, \Gamma(\ell+2\xi)}\, \int_0^\infty\, \frac{J_{2q-1}(2\sqrt{st})}{(st)^{q-1/2}}\, \frac{s^\ell}{\ell!}\, dm_q(s).
\]
Then, a classical integral equality involving Laguerre polynomials (see \cite[vol. II, p. 190]{E}) implies that
\[
\int_0^\infty\, \frac{J_{2q-1}(2\sqrt{st})}{(st)^{q-1/2}}\, \frac{s^\ell}{\ell!}\, dm_q(s) = e^{-t}\, L_\ell^{2q-1}(t).
\]
Hence we obtain
\[
\Big( N_q(L_n^{2\xi-1})\, ,\, L_k^{2\xi-1} \Big)_\xi = \sum_{\ell=0}^n\, \frac{(-1)^\ell\, \Gamma(n+2\xi)}{(n-\ell)!\, \Gamma(\ell+2\xi)}\, \Big( e^{-t}\, L_\ell^{2q-1}(t)\, ,\, L_k^{2\xi-1}(t) \Big)_\xi.
\]
Finally, we apply \cite[eq. (14)]{conti} to get
\[
\begin{aligned}
 & \Big( e^{-t}\, L_\ell^{2q-1}(t)\, ,\, L_k^{2\xi-1}(t) \Big)_\xi = \int_0^\infty\, L_\ell^{2q-1}(t)\, L_k^{2\xi-1}(t)\, \chi_{2\xi-1}(t)\, e^{-2t}\, dt = \\
 & = \frac{\Gamma(2\xi)}{2^{2\xi}}\, \sum_{j=0}^{\min\{\ell,k\}}\, \frac{\Gamma(\ell+2q)}{(\ell-j)! \, \Gamma(j+2q)}\, \frac{\Gamma(k+2\xi)}{(k-j)! \, \Gamma(j+2\xi)}\, \frac{\Gamma(j+2\xi)}{j! \, \Gamma(2\xi)}\, 2^{-2j} \cdot \\
 & \cdot\ \hyper\left( -\ell+j,2\xi +j;2q+j;\frac 12\right)\, \hyper\left( -k+j,2\xi +j;2\xi+j;\frac 12\right) =\\
 & = \sum_{j=0}^{\min\{\ell,k\}}\, \frac{2^{-(k+2\xi+j)} \Gamma(\ell+2q) \Gamma(k+2\xi) }{j!\, (\ell-j)! \, (k-j)! \, \Gamma(j+2q)}\, \hyper\left( -\ell+j,2\xi +j;2q+j;\frac 12\right),
\end{aligned}
\]
where in the last equation we have used that $\hyper(\alpha,\beta;\beta;z) = (1-z)^{-\alpha}$.
\end{proof}

Using Lemmas \ref{lem:m-matr} and \ref{lem:n-matr} in \eqref{eq-prob-0} and \eqref{inf-mat-0-impl}, we obtain
\begin{theorem}\label{thm:mat-0}
Assume $\xi:=\Re(q)>0$ and $q\not= 1/2$. Then equation \eqref{the-problem-0} has a nonzero solution $\phi \in L^2(m_\xi)$ if and only if there exists a nonzero infinite vector $\Phi := (\Phi_n)$ such that
\begin{equation}\label{sistema}
\left\{ \begin{array}{l}
\sum_{n=0}^\infty\, |\Phi_n|^2\, \frac{\Gamma(n+2\xi)}{n!}<\infty \\[0.2cm]
A^\pm_q \Phi = D_q \Phi
\end{array} \right.
\end{equation}
where $A^\pm_q = (a^\pm_{kn}(q))$ and $D_q = (d_{kn}(q))$ are infinite matrices with $k,n\in \N_0$, given by
\[
\begin{aligned}
a^\pm_{kn}(q) =\, & \frac{\Gamma(k+n+2\xi)}{k!\, n!}\, 2^{-(k+n+2\xi)} \pm \frac{\Gamma(n+2\xi)\Gamma(k+2\xi)}{2^{k+2\xi}}\, \sum_{\ell=0}^n\,  \frac{(-1)^\ell\, \Gamma(\ell+2q)}{(n-\ell)!\, \Gamma(\ell+2\xi)}\, \cdot \\
& \cdot \left(\sum_{j=0}^{\min\{\ell,k\}}\, \frac{2^{-j}}{j!\, (\ell-j)! \, (k-j)! \, \Gamma(j+2q)}\, \hyper\left( -\ell+j,2\xi +j;2q+j;\frac 12\right)\right),\\[0.3cm]
d_{kn}(q) =\, & \frac{\Gamma(k+2\xi)}{k!}\, \delta_{k,n}.
\end{aligned}
\]
If the solution $\Phi := (\Phi_n)$ to system \eqref{sistema} exists, then the solution to \eqref{the-problem-0} can be written as
\[
\phi(t) = \sum_{n=0}^\infty\, \Phi_n\, L_n^{2\xi-1}(t),
\]
and viceversa.
\end{theorem}

We now prove the analogous result for \eqref{the-problem-b}. By Proposition \ref{prop-m-nq} we have $N_q(\chi_{-1}) \in L^2(m_\xi)$, and using the definition of the operator $M$ one can show that $(I-M)\chi_{-1} \in L^2(m_\xi)$. Therefore $(I-M-N_q)\chi_{-1}  \in L^2(m_\xi)$ and, using the basis of Laguerre polynomials, we can write that a function $\phi\in L^2(m_\xi)$ is a nonzero solution of 
\eqref{the-problem-b} if and only if
\[
\Big( (M+ N_q - I)\phi\, ,\, L^{2\xi-1}_k\Big)_\xi = \Big((I-M-N_q)\chi_{-1}\, ,\, L_k^{2\xi-1}\Big)_\xi, \quad \forall\, k\in \N_0.
\]
Writing $\phi$ as in \eqref{expansion} we obtain again a linear system for an infinite vector $\Phi=(\Phi_n)$ which in this case is not homogeneous. Let's write the expansion \eqref{expansion} for $(I-M-N_q)\chi_{-1}$ as
\begin{equation}\label{expansion-term}
((I-M-N_q)\chi_{-1})(t) = \sum_{n=0}^\infty\, \Psi_n\, L_n^{2\xi-1}(t),
\end{equation}
where $\{\Psi_n\}_n$ is the sequence in $\C$ which satisfies
\begin{equation}\label{condition-term}
\begin{aligned}
& \Big((I-M-N_q)\chi_{-1}\, ,\, L_n^{2\xi-1}\Big)_\xi = \Psi_n \frac{\Gamma(n+2\xi)}{n!},\\
& \Big\|(I-M-N_q)\chi_{-1}\Big\|_2^2 =  \sum_{n=0}^\infty\, |\Psi_n|^2\, \frac{\Gamma(n+2\xi)}{n!} <+\infty.
\end{aligned}
\end{equation}
Then, if $A_q^+$ and $D_q$ are the infinite matrices defined in \eqref{inf-mat-0-impl} and \eqref{inf-mat-diag}, there exists a nonzero solution to \eqref{the-problem-b} if and only if there exists a nonzero infinite vector $\Phi := (\Phi_n)$ such that
\begin{equation}\label{eq-prob-b}
\left\{ \begin{array}{l}
\sum_{n=0}^\infty\, |\Phi_n|^2\, \frac{\Gamma(n+2\xi)}{n!}<\infty \\[0.2cm]
A^+_q \Phi = D_q \Phi +D_q \Psi
\end{array} \right.
\end{equation}
where $\Psi :=(\Psi_n)$ is the infinite vector defined in \eqref{expansion-term} and \eqref{condition-term}. The coefficients of $A_q^+$ are given in Theorem \ref{thm:mat-0}, hence to have an explicit expression for the linear system \eqref{eq-prob-b} we only need to have an expression for the $\Psi_n$'s. This is the content of our last result.

\begin{theorem}\label{thm:mat-b}
Assume $\xi:=\Re(q)>0$ and $q\not= 1/2$. Then equation \eqref{the-problem-b} has a nonzero solution $\phi \in L^2(m_\xi)$ if and only if there exists a nonzero infinite vector $\Phi := (\Phi_n)$ satisfying \eqref{eq-prob-b} where $A^+_q = (a^+_{kn}(q))$ is as in Theorem \ref{thm:mat-0}, $D_q = (d_{kn}(q))$ is defined in \eqref{inf-mat-diag}, and $\Psi=(\Psi_n)$ is the infinite vector given by
\[
\Psi_n = \int_{1/2}^1\, t^{2\xi-2} (1-t)^n\, dt \ - \frac{1}{2q-1}\, \int_0^1\, \frac{t^{n+2q-2}}{(1+t)^{n+2\xi+1}}\, \Big( (2\xi+4q)\, t -n - 2(2\xi+1)\, \frac{t^2}{1+t} + 2n\, \frac{t}{1+t}\Big)\, dt.
\]
If the solution $\Phi := (\Phi_n)$ to system \eqref{eq-prob-b} exists, then the solution to \eqref{the-problem-b} can be written as
\[
\phi(t) = \sum_{n=0}^\infty\, \Phi_n\, L_n^{2\xi-1}(t),
\]
and viceversa.
\end{theorem}

\begin{proof}
We only need to compute the terms
\[
\Psi_n = \frac{n!}{\Gamma(n+2\xi)}\, \Big((I-M-N_q)\chi_{-1}\, ,\, L_n^{2\xi-1}\Big)_\xi.
\]

We begin with the first part. Using \eqref{laguerre} we can write
\[
\begin{aligned}
& \Big((I-M)\chi_{-1}\, ,\, L_n^{2\xi-1}\Big)_\xi = \sum_{\ell=0}^n\, \frac{\Gamma(n+2\xi)}{(n-\ell)!\, \Gamma(\ell+2\xi)}\, \frac{(-1)^\ell}{\ell!}\, \int_0^\infty\, \frac{1-e^{-t}}{t}\, t^\ell\, dm_\xi(t)\\
& = \sum_{\ell=0}^n\, \frac{\Gamma(n+2\xi)}{(n-\ell)!\, \Gamma(\ell+2\xi)}\, \frac{(-1)^\ell}{\ell!}\, \Big( \int_0^\infty\, e^{-t}\, t^{2\xi +\ell-2}\, dt - \int_0^\infty\, e^{-2t}\, t^{2\xi +\ell-2}\, dt \Big) \\
& = \sum_{\ell=0}^n\, \frac{\Gamma(n+2\xi)}{(n-\ell)!\, \Gamma(\ell+2\xi)}\, \frac{(-1)^\ell}{\ell!}\, \Big( \Gamma(\ell + 2\xi -1) - \frac{\Gamma(\ell + 2\xi -1)}{2^{\ell + 2\xi -1}} \Big),
\end{aligned}
\]
where we have used the definition of the Gamma function and a change of variables in the second integral. Then we use
\[
\frac{\Gamma(\ell + 2\xi -1)}{\Gamma(\ell + 2\xi)} = \frac{1}{\ell + 2\xi -1} \quad \text{and} \quad t^{2\xi-2} (1-t)^n = \sum_{\ell=0}^n\, {n\choose \ell}\, (-1)^\ell\, s^\ell,
\]
to obtain
\begin{equation} \label{1/t-prima}
\begin{aligned}
\Big((I-M)\chi_{-1}\, ,\, L_n^{2\xi-1}\Big)_\xi & = \frac{\Gamma(n+2\xi)}{n!}\, \sum_{\ell=0}^n\, {n\choose \ell}\, \frac{(-1)^\ell}{\ell+2\xi -1}\, \left( 1- \frac{1}{2^{\ell + 2\xi -1}} \right)\\
& = \frac{\Gamma(n+2\xi)}{n!}\, \int_{1/2}^1\, t^{2\xi-2} (1-t)^n\, dt.
\end{aligned}
\end{equation}

Then we compute an explicit expression for the term containing $N_q(\chi_{-1})$. For these computations we refer to formulas in \cite{nist}. First of all we use the definition of the action of $N_q$ in \eqref{eq-m-nq} and \cite[Eq. (8.6.2)]{nist} to write
\[
N_q(\chi_{-1})(t) = t^{1-2q}\, \gamma(2q-1,t),
\]
where $\gamma(\cdot,\cdot)$ denotes the \emph{incomplete Gamma function} (see e.g. \cite[Chap. 8]{nist}). We can then use \eqref{laguerre} to write
\begin{equation}\label{1/t-seconda-prima}
\begin{aligned}
\Big(N_q(\chi_{-1})\, ,\, L_n^{2\xi-1}\Big)_\xi & = \sum_{\ell=0}^n\, \frac{\Gamma(n+2\xi)}{(n-\ell)!\, \Gamma(\ell+2\xi)}\, \frac{(-1)^\ell}{\ell!}\, \int_0^\infty\, t^{\ell+1-2q}\, \gamma(2q-1,t)\, dm_\xi(t) =\\
& = \sum_{\ell=0}^n\, \frac{\Gamma(n+2\xi)}{(n-\ell)!\, \Gamma(\ell+2\xi)}\, \frac{(-1)^\ell}{\ell!}\, \int_0^\infty\, t^{\ell+2(\xi-q)}\, \gamma(2q-1,t)\, e^{-t}\, dt =\\
& = \frac{1}{2q-1}\, \frac{\Gamma(n+2\xi)}{n!}\, \sum_{\ell=0}^n\, {n\choose \ell}\, \frac{(-1)^\ell}{2^{\ell+2\xi}}\, \hyper\left(\ell+2\xi,1;2q;\frac 12 \right),
\end{aligned}
\end{equation}
where in the last line we have used the integral equality in \cite[Eq. (8.14.5)]{nist}. At this point we use some properties of the hypergeometric function to obtain a formulation in terms of the integral of elementary functions, just as in the term containing $(I-M)\chi_{-1}$. In particular we use \cite[Eq. (15.8.1) and (15.5.19)]{nist} which imply
\[
\begin{aligned}
& \frac{1}{2^{\ell+2\xi}}\, \hyper\left(\ell+2\xi,1;2q;\frac 12 \right) =\, \hyper\left(\ell+2\xi,2q-1;2q;-1\right) \\
& = -\frac{2(\ell+2\xi+1)}{2q+1}\, \hyper\left(\ell+2\xi+2,2q+1;2q+2;-1\right) + \frac{\ell+2\xi+4q}{2q}\, \hyper\left(\ell+2\xi+1,2q;2q+1;-1\right).
\end{aligned}
\]
So that we can apply the integral representation \cite[Eq. (15.6.1)]{nist} given by
\[
\hyper(a,b;c;z)= \frac{\Gamma(c)}{\Gamma(b)\, \Gamma(c-b)}\, \int_0^1\, t^{b-1}\, (1-t)^{c-b-1}\, (1-zt)^{-a}\, dt
\]
for $z\in \C \setminus [1,\infty)$ and $\Re(c)>\Re(b)>0$. Therefore
\[
\begin{aligned}
\frac{1}{2^{\ell+2\xi}}\, \hyper\left(\ell+2\xi,1;2q;\frac 12 \right) = & -\frac{2(\ell+2\xi+1)}{2q+1}\, \frac{\Gamma(2q+2)}{\Gamma(2q+1)\, \Gamma(1)}\, \int_0^1\, t^{2q}\, (1+t)^{-\ell-2\xi-2}\, dt \\ & + \frac{\ell+2\xi+4q}{2q}\, \frac{\Gamma(2q+1)}{\Gamma(2q)\, \Gamma(1)}\, \int_0^1\, t^{2q-1}\, (1+t)^{-\ell-2\xi-1}\, dt \\ = & -2(\ell+2\xi+1)\, \int_0^1\, t^{2q}\, (1+t)^{-\ell-2\xi-2}\, dt \\ & + (\ell+2\xi+4q)\, \int_0^1\, t^{2q-1}\, (1+t)^{-\ell-2\xi-1}\, dt. 
\end{aligned}
\]
Using this expression in \eqref{1/t-seconda-prima} along with the identities
\[
\begin{aligned}
& \sum_{\ell=0}^n\, {n\choose \ell}\, (-1)^\ell\, (\ell+2\xi+1)\, (1+t)^{-\ell-2\xi-2} = \frac{-n\, t^{n-1}+ (2\xi+1)\, t^n}{(1+t)^{n+2\xi+2}}, \\
& \sum_{\ell=0}^n\, {n\choose \ell}\, (-1)^\ell\, (\ell+2\xi+4q)\, (1+t)^{-\ell-2\xi-1} = \frac{-n\, t^{n-1}+ (2\xi+4q)\, t^n}{(1+t)^{n+2\xi+1}},
\end{aligned}
\]
gives
\begin{equation}\label{1/t-seconda}
\begin{aligned}
& \Big(N_q(\chi_{-1})\, ,\, L_n^{2\xi-1}\Big)_\xi \\
& =  \frac{1}{2q-1}\, \frac{\Gamma(n+2\xi)}{n!}\, \int_0^1\, \Big( \frac{2n\, t^{n+2q-1}- 2(2\xi+1)\, t^{n+2q}}{(1+t)^{n+2\xi+2}} + \frac{-n\, t^{n+2q-2}+ (2\xi+4q)\, t^{n+2q-1}}{(1+t)^{n+2\xi+1}} \Big)\, dt\\
& = \frac{1}{2q-1}\, \frac{\Gamma(n+2\xi)}{n!}\, \int_0^1\, \frac{t^{n+2q-2}}{(1+t)^{n+2\xi+1}}\, \Big( (2\xi+4q)\, t -n - 2(2\xi+1)\, \frac{t^2}{1+t} + 2n\, \frac{t}{1+t}\Big)\, dt.
\end{aligned}
\end{equation}
The proof is finished by putting together \eqref{1/t-prima} and \eqref{1/t-seconda}.
\end{proof}

We believe that to have explicit expressions for the matrices and the vectors in \eqref{eq-prob-0} and \eqref{eq-prob-b} will turn useful for a numerical approach to the two linear systems.


\end{document}